\documentclass{article}
\usepackage{amssymb,amsmath,mathrsfs,latexsym,amscd}
\usepackage{exscale,latexsym,amsthm,graphics}

\usepackage{makeidx}
\makeindex
\usepackage{epsfig}
\usepackage{txfonts}
\newtheorem{thm}{Theorem}[section]
\newtheorem{cor}[thm]{Corollary}

\newtheorem{lemma}[thm]{Lemma}
\newtheorem{prop}[thm]{Proposition}

\newtheorem{definition}[thm]{Definition}

\newtheorem{remark}[thm]{Remark}

\def\bM {{\mathbb M}}

\def\qed{{\hfill $\Box$ \bigskip}}

\def\R {{\mathbb R}}

\def\EE{{\mathbb E}}
\def\P{{\mathbb P}}
\newcommand{\F}{\mathcal{F}}

\def\E{{\mathcal E}}
\def\P{{\mathbb P}}

\numberwithin{equation}{section}
\begin{document}

\noindent
{{\Large\bf Liouville distorted Brownian motion}}

\bigskip
\noindent
{\bf Jiyong Shin}
\\

\noindent
{\small{\bf Abstract.} The Liouville Brownian motion  was introduced in \cite{GRV} as a time changed process $B_{A_t^{-1}}$ of a planar Brownian motion $(B_t)_{t \ge 0}$, where $(A_t)_{t \ge 0}$ is the positive continuous additive functional  of  $(B_t)_{t \ge 0}$  in the strict sense w.r.t. the Liouville measure. We first consider a distorted Brownian motion $(X_t)_{t\ge0}$ starting from all points in $\R^2$ associated to a Dirichlet form $(\E, D(\E))$ (see \cite{ShTr14}). We show that  the positive continuous additive functional  $(F_t)_{t \ge 0}$ of  $(X_t)_{t \ge 0}$  in the strict sense w.r.t. the Liouville distorted measure can be constructed.\\



\section{Introduction}
The Liouville Brownian motion, introduced by C. Garban, R. Rhodes, V. Vargas in \cite{GRV}, is a Markov process defined as a time changed process on $\R^2$. 
 By classical theory of Gaussian multiplicative chaos (cf. \cite{Kah}), the Liouville measure $M_{\gamma}$, $\gamma \in (0,2)$ is well defined (see Section \ref{ss;mlgm}). In \cite{GRV} the positive continuous additive functional $(A_t)_{t \ge 0}$ of a planar Brownian motion $(B_t)_{t \ge 0}$ in the strict sense w.r.t. the Liouville measure $M_{\gamma}$ is constructed and then the Liouville Brownian motion is defined as $B_{A_t^{-1}}$. 

In this paper we are concerned with the extension of the Liouville Brownian motion $B_{A_t^{-1}}$ to more general Markov processes. Note that the planar Brownian motion $(B_t)_{t \ge 0}$ is associated with the Dirichlet form
\[
\E^{\prime}(f,g):= \frac{1}{2} \int_{\R^2}   \langle \nabla f,  \nabla g \rangle \ dx, 
\]
where $f, g \in \mathcal{F} = \{f \in L^2(\R^2,dx) \mid \nabla f \in L^2(\R^2,dx) \}$. As an extension of the Dirichlet form $(\E^{\prime},\mathcal{F})$, we first consider a more general Dirichlet form   $(\E,D(\E))$  which is defined as the closure of  the symmetric bilinear form
\[
\E(f,g):= \frac{1}{2} \int_{\R^2}   \langle \nabla f(x),  \nabla g(x) \rangle \ \rho(x) \ dx, \quad f, g \in C_0^{\infty}(\R^2) ,
\]
on $L^2(\R^2,\rho dx)$ where   $\rho(x ):= |x|^{\alpha}$, $\alpha \in (-2, \infty)$ (see \cite{ShTr14}).
It is known from  \cite{ShTr13}  and \cite[Section 3]{ShTr14} that there exists the distorted Brownian motion $( (X_t)_{t \ge 0}, (\P_x)_{x \in \R^2} )$ starting from all points in $\R^2$ associated with the Dirichlet form $(\E,D(\E))$. Then using the estimates of the resolvent kernel and part Dirichlet form method,  we can construct the positive continuous additive functional $(F_t)_{t \ge 0}$ of $(X_t)_{t \ge 0}$  in the strict sense  w.r.t.  the Liouville distorted measure $M_{\gamma}^{\rho}$ (see Section \ref{ss;mlgm}). Similarly to the Liouville Brownian motion, the Liouville distorted Brownian motion is defined as $X_{F_t^{-1}}$.\\

\noindent{\bf{Notations}}:\\
For any open set $U$ in $\R^2$, we denote the set of all  Borel measurable functions and the set of all bounded Borel measurable functions on $U$  by  $\mathcal{B}(U)$ and $\mathcal{B}_b(U)$, respectively. The usual $L^q$-spaces $L^q(U, m)$, $q \in[1,\infty]$ are equipped with $L^{q}$-norm $\| \cdot \|_{L^q (U,m) }$ with respect to the measure $m$ on $U$ and $\mathcal{A}_b$ : = $\mathcal{A} \cap \mathcal{B}_b(U)$ for $\mathcal{A} \subset L^q(U,m)$.  The inner product on $L^2(U, m)$ is denoted by $(\cdot,\cdot)_{L^2(U, m)}$.  The indicator function of a set $A$ is denoted by $1_A$. Let $\nabla f : = ( \partial_{1} f, \dots , \partial_{d} f )$  and  $\Delta f : = \sum_{j=1}^{d} \partial_{jj} f$ where $\partial_j f$ is the $j$-th weak partial derivative of $f$ and $\partial_{jj} f := \partial_{j}(\partial_{j} f) $, $j=1, \dots, d$.  As usual $dx$ denotes the Lebesgue measure on $\R^2$. Here $C_0^{\infty}(\R^2)$ denotes the set of all infinitely differentiable functions with compact support in $\R^2$. We equip $\R^2$ with the Euclidean norm $|\cdot|$ and the corresponding inner product $\langle \cdot, \cdot \rangle$.

\section{Massive Gaussian free field and Gaussian multiplicative chaos}\label{ss;mlgm}
We first state the definition of the massive Gaussian free field as stated in \cite{GRV}. The massive Gaussian free field on $\R^2$ is a centered Gaussian random distribution (in the sense of Schwartz)  on a probability space $(\Omega, \mathcal{A}, P)$ with covariance function given by the Green function $G^{(m)}$ of the operator $m^2 -\Delta$, $m > 0$, i.e.
\[
(m^2 - \Delta) G^{(m)}(x, \cdot) = 2 \pi \delta_x, \quad x \in \R^2, 
\] 
where $\delta_x$ stands for the Dirac mass at $x$.
The massive Green function with the operator $(m^2 - \Delta)$ can be written as 
\[
G^{(m)}(x,y) = \int_0^{\infty} e^{- \frac{m^2}{2}s - \frac{\| x-y \|^2}{2s}} \frac{ds}{2s} 
= \int_{1}^{\infty} \frac{k_{m} (s (x-y)) }{s} \ ds, \quad x,y \in \R^2,
\]
where 
\[
k_m (z) = \frac{1}{2} \int_0^{\infty} e^{- \frac{m^2}{2s} \| z  \|^2 - \frac{s}{2} } \ ds.
\]

Let $(c_n)_{n \ge 1}$ be an unbounded strictly increasing sequence such that $c_1 = 1$ and  $(Y_n)_{n \ge 1}$ be a family of  independent centered continuous Gaussian fields on $\R^2$ on the probability space $(\Omega, \mathcal{A}, P)$ with covariance kernel given by 
\[
E[Y_n (x) \ Y_n (y)] = \int_{c_{n-1}}^{c_n}  \frac{k_m (s(x-y))}{s} ds.
\]
The massive Gaussian free field is the Gaussian distribution defined by 
\[
X(x) = \sum_{k \ge 1} Y_k (x).
\]

We define $n$-regularized field by
\[
X_n(x) = \sum_{k=1}^{n} Y_k(x), \quad n \ge 1
\]
and the associated  $n$-regularized measure by
\[
M_{n,\gamma}^{\rho} (dz)=  \exp \Big( \gamma X_n(z) - \frac{\gamma^2}{2} E[X_n(z)^2] \Big) \ \rho(z) dz, \quad \gamma \in (0,2),
\]
where $\rho dz$ is a positive Radon measure on $\R^2$. By the classical theory of Gaussian multiplicative chaos (see \cite{Kah}), $P$-a.s. the family $(M_{n,\gamma})_{n \ge 1}$ weakly converges to the Liouville distorted measure
\[
M_{\gamma}^{\rho} (dz) =  \exp \Big( \gamma X(z) - \frac{\gamma^2}{2} E[X(z)^2] \Big) \ \rho(z) dz.
\]
It is known from \cite{Kah} that $M_{\gamma}^{\rho}$ is a Radon measure on $\R^2$. 
If $\rho =1$, we denote the n-regularized Liouville measure $M_{n,\gamma}^{\rho}$ and the Liouville measure $M_{\gamma}^{\rho}$ by $M_{n,\gamma}$ and $M_{\gamma}$, respectively.
\section{Liouville distorted Brownian motion}

We consider  $\rho(x ):= |x|^{\alpha}$, $\alpha \in (-2, \infty)$ and the symmetric bilinear form
\[
\E(f,g):= \frac{1}{2} \int_{\R^2}   \langle \nabla f(x),  \nabla g(x) \rangle \ m (dx), \quad f, g \in C_0^{\infty}(\R^2) ,
\]
where $m:= \rho dx$. It is known that $(\E,C_0^{\infty}(\R^2) )$ is closable and its closure is a strongly local, regular Dirichlet form $(\E,D(\E))$ on $L^2(\R^d,m)$ (cf. e.g. \cite{ShTr14}). Let  $(T_t)_{t > 0}$ and $(G_{\beta})_{\beta > 0}$ be the $L^2(\R^2, m)$-semigroup and resolvent associated to $(\E,D(\E))$ (see \cite{FOT}). By \cite{ShTr13}  and \cite[Section 3]{ShTr14} there exists the distorted Brownian motion associated with the Dirichlet form $(\E,D(\E))$
\[
\bM : = (\Omega, \F, (\F_t)_{t \ge 0}, \zeta, (X_t)_{t \ge 0}, (\P_x)_{x \in \R^2} )
\]
with transition function $(P_t)_{t \ge 0}$ where  $\zeta$ is the lifetime.
Moreover, it is known from \cite[Section 3]{ShTr14} that  there exists a jointly continuous transition kernel density  $p_{t}(x,y)$ such that
\[
P_t f(x) := \int_{\R^2} p_t(x,y)\, f(y) \ m(dy), \quad t>0, \ x,y \in \R^2, \ f \in \mathcal{B}_b(\R^2)
\]
is an $m$-version of $T_t f$ if $f  \in  L^2(\R^2 , m)_b$. We set $P_0 : = id$.
Taking the Laplace transform of $p_{\cdot}(x, y)$, we obtain a $\mathcal{B}(\R^2) \times \mathcal{B}(\R^2)$ measurable non-negative resolvent kernel density $r_{\beta}(x,y)$ such that
\[
R_{\beta} f(x) := \int_{\R^2} r_{\beta}(x,y)\, f(y) \, m(dy), \quad \beta>0, \ x \in \R^2, f \in \mathcal{B}_b(\R^2),
\]
is an $m$-version of $G_{\beta} f$ if $f  \in  L^2(\R^2, m)_b$.  

We present some definitions and properties concerning  $(\E,D(\E))$. We will refer to \cite{FOT} till the end, hence some of its standard notations may be adopted below without definition. For any set $A \subset \R^2$ the capacity of $A$ is defined as
\[
\text{Cap(A)} = \inf_{\begin{subarray} \ \text{open} \ B \subset \R^2  \\    A \subset B \end{subarray}} \inf_{\begin{subarray} \  f \in D(\E) \\ 1_{B} \cdot f \ge 1 \ \text{m-a.e.}  \end{subarray}} \E(f,f) +  (f,f)_{L^2(\R^2,m)}.
\]
\begin{definition}
Let $B$ be an open set in $\R^2$. For $x \in B, t \ge 0, \beta>0$ and $p \in [1,\infty)$  let
\begin{itemize}
\item $\sigma_{B^c} := \inf\{t>0 \,|\;X_t \in B^c\},   \quad   D_{B^c} := \inf\{t\ge0 \,|\;X_t \in B^c\}$,
\item $P^{B}_{t}f(x) : = \EE_x [f(X_t) ; t<\sigma_{B^c} ], \quad f \in \mathcal{B}_{b}(B)$,
\item $R^{B}_{\beta}f(x) : = \EE_x \big[\int_{0}^{\sigma_{B^c}} e^{-\beta s} f(X_s) \, ds \big],  \quad  f \in \mathcal{B}_{b}(B)$ ,
\item $D(\E^{B}): = \{u \in D(\E) \, | \; u=0 \,\,  \mathcal{E}\text{-q.e} \; on  \; B^c \}$.
\item $\E^{B} : = \mathcal{E} \, |_{D(\E^{B})\times D(\E^{B})}$.
\item $\E^{B}_1(f,g) : = \E^{B}(f,g) + \int_{B} f g \ dm,  \quad     f,g \in D(\E^{B})$.
\item $\| \,f\, \|_{D(\mathcal{E}^{B})} : = \E^{B}_1(f,f)^{1/2},  \quad      f \in D(\E^{B})$.
\end{itemize}
\end{definition}
It is known that  $(\E^{B},D(\E^{B}))$  is a regular Dirichlet form on $L^2(B, m)$, which is called the part Dirichlet form of $(\E, D(\E))$ on $B$ (cf. \cite[Section 4.4]{FOT}). Let  $(T^{B}_t)_{t > 0}$ and $(G^{B}_{\beta})_{\beta > 0}$ be the $L^2(B, m)$-semigroup and resolvent associated to $(\E^{B},D(\E^{B}))$.
Then $P^{B}_{t} f, \;  R^{B}_{\beta}f$ is an $m$-version of $T^{B}_t f,  G^{B}_{\beta}f$, respectively for any $f \in L^2(B,m)_b$. Since $P_t^{B} 1_{A}(x) \le P_t 1_{A}(x)$ for any $ A\in \mathcal{B}(B)$, $x \in B$ and $m$ has full support on $E$, $A \mapsto P_t^{B} 1_{A}(x), \; A \in \mathcal{B}(B)$ is absolutely continuous with respect to $1_{B} \cdot m$. Hence there exists a (measurable) transition kernel density $p^{B}_{t}(x, y)$, $x,y \in B$, such that
\begin{equation}\label{pdldbd}
P_t^{B} f(x) = \int_{B} p_t^{B} (x,y) \,  f(y) \, m(dy) ,\; t>0 \;, \;\; x\in B
\end{equation}
for $f  \in  \mathcal{B}_{b}(B)$. Correspondingly, there exists a (measurable) resolvent kernel density $r_{\beta}^{B}(x,y)$, such that
\[
R_{\beta}^{B} f(x) = \int_{B} r_{\beta}^{B} (x,y) \  f(y) \ m(dy), \quad \beta>0, \quad x \in B
\]
for $f  \in \mathcal{B}_{b}(B)$.
For a signed Radon measure $\mu$ on $B$, let us define
\[
R_{\beta}^{B} \mu (x) = \int_{B} r_{\beta}^{B} (x,y) \ \mu(dy), \quad \beta>0, \quad x \in B
\]
whenever this makes sense.
The process defined by
\[
X^{B}_t(\omega)=
\begin{cases}
X_t(\omega), \quad 0\le t < D_{B^c} (\omega) \\
\Delta, \quad t \ge D_{B^c} (\omega)
\end{cases}
\]
is called the part process associated to $\E^{B}$ and is denoted by $\bM|_{B}$.  The part process $\bM|_{B}$ is a Hunt process on $B$ (see \cite[p.174 and Theorem A.2.10]{FOT}).  In particular, by \eqref{pdldbd}
$\bM|_{B}$ satisfies the absolute continuity condition on $B$. 

A positive Radon measure $\mu$ on $B$ is said to be of finite energy integral if
\[
\int_{B} |f(x)|\, \mu (dx) \leq c \sqrt{\E^{B}_1(f,f)}, \; f\in D(\E^{B}) \cap C_0(B),
\]
where $c$ is some constant independent of $f$ and $C_0(B)$ is the set of all compactly supported continuous functions on $B$. A positive Radon measure $\mu$ on $B$ is of finite energy integral (on $B$) if and only if there exists a unique function $U_{1}^{B} \, \mu\in D(\E^{B} )$ such that
\[
\E^{B}_{1}(U_{1}^{B} \, \mu, f) = \int_{B} f(x) \, \mu(dx),
\]
for all $f \in D(\E^{B}) \cap C_0(B)$. $U_{1}^{B} \, \mu$ is called $1$-potential of $\mu$. In particular, $R_{1}^{B} \mu$ is a version of $U_{1}^{B} \mu$ (see e.g. \cite[Exercise 4.2.2]{FOT}). The measures of finite energy integral are denoted by $S_0^{B}$. We further define $S_{00}^{B} : = \{\mu\in S_0^{B}  \mid  \mu(B)<\infty, \|U_{1}^{B} \mu\|_{L^{\infty}(B,m)}<\infty \}$.

If  $\mu \in S_{00}^{B}$,
then there exists a unique $A \in A_{c,1}^{+, B}$ with $\mu = \mu_{A}$, i.e. $\mu$ is the Revuz measure of $A$ (see \cite[Theorem 5.1.6]{FOT}). Here, $A_{c,1}^{+,B}$ denotes the positive continuous additive functionals on $B$ in the strict sense.\\

We define
\[
E_k:= \{x \in \R^2  \mid 1 /k <  |x| <k \}, \quad \text{and} \quad  E:= \bigcup_{k \ge 1} E_k  = \R^2 \setminus \{0\}.
\]
The following proposition recalls some properties of the Dirichlet form $(\E,D(\E))$ as stated in \cite[Theorem 2.10, Lemma 3.13]{ShTr14}:
\begin{prop}\label{p;cacod}
\begin{itemize}
\item[(i)]  For any $\alpha \in (-2, \infty)$, the Dirichlet form $(\E,D(\E))$ is conservative.
\item[(ii)]  For any $\alpha \in [0,\infty)$, $\emph{Cap}(\{0\})=0$.
\item[(iii)] For any $\alpha \in [0,\infty)$ and  any $x \in E$
\[
\P_x \Big(\lim_{k \rightarrow \infty} D_{E_k^c} = \infty \Big) = \P_x \Big(\lim_{k \rightarrow \infty} \sigma_{E_k^c} = \infty \Big)=1.
\]
\end{itemize}
\end{prop}

From now on till the end of this paper, we consider
\[
\rho (x ) = |x|^{\alpha}, \quad \alpha \in [2, \infty).
\]
\begin{lemma}\label{l:srh} 
Almost surely in $X$, for any relatively compact open set $G  \subset E$,
\[
M_{\gamma}^{\rho} (G) \le \sup_{x \in G} \rho(x) \ M_{\gamma} (G).
\]
\end{lemma}
\proof
The statement follows from
\[
M_{n,\gamma}^{\rho} (G) \le \sup_{x \in G} \rho(x)  \int_{G} \exp \Big( \gamma X_n(z) - \frac{\gamma^2}{2} E[X_n(z)^2] \Big)  \ dz.
\]
\qed

\begin{lemma}\label{l;rdldb} 
Let $G$ be any relatively compact open set in $E$ with $G\subset \overline{G} \subset E$. For any $x,y \in G$ and any $\delta >0$
\[
r_1^{G} (x,y) \le \  c_1 \frac{1}{|x-y|^{\delta}},
\]
where $c_1>0$ is some constant.
\end{lemma}
\proof
By \cite[Lemma 3.10]{ShTr14} for $m$-a.e. $x,y \in G$ and any $\delta >0$
\[
r_1^{G} (x,y) \le \  c_1 \frac{1}{|x-y|^{\delta}},
\]
where $c_1>0$ is some constant. Since $p_t(x,y)$ is jointly continuous, the statement holds for all $x,y \in G$.
\qed
\begin{thm}\label{t;s0ds}
Almost surely in $X$, for any relatively compact open set $G\subset \overline{G} \subset E$, $1_{G} \cdot M^{\rho}_{\gamma} \in S^G_{00}$.
\end{thm}
\proof
Clearly, $M^{\rho}_{\gamma} (G) < \infty$. By Lemma \ref{l:srh} and Lemma \ref{l;rdldb}, for any $x \in G$ and $\delta > 0$,
\begin{eqnarray*}
R_1^G (1_{G} \cdot M^{\rho}_{\gamma}) (x) &\le&  c_1  \int_{G} \frac{1}{|x-y|^{\delta}} \ M_{\gamma}^{\rho}(dy) \\
&\le& c_1\ \sup_{y \in G} \rho(y)  \int_{G} \frac{1}{|x-y|^{\delta}} \ M_{\gamma}(dy).
\end{eqnarray*}
Since  $G$ is a relatively compact open set, we can find a constant $R>0$ such that $G \subset B_R(0)$.
 By \cite[Theorem 2.2]{GRV}, there exist a constant $c_2>0$ and $\alpha>0$ (depending on $R>0$ and $\gamma$) such that
for all $r \in (0,R)$ and $x \in B_R(0)$
\[
M_{\gamma}(B_r(x)) \le c_2 \ r^{\alpha}.
\]
By taking $0 < \delta < \alpha$, we obtain
\begin{eqnarray*}
 \int_{G} \frac{1}{|x-y|^{\delta}} \ M_{\gamma}(dy) &\le& \int_{B_R(0)} \frac{1}{|x-y|^{\delta}} \ M_{\gamma}(dy)  \le \sum_{n \ge 0} \int_{B_R(0) \cap \{2^{-n} R < |x-y| \le 2^{-n+1} R \} } \frac{1}{|x-y|^{\delta}} \ M_{\gamma}(dy) \\
&\le&  \sum_{n \ge 0} 2^{ n\delta} R^{-\delta} M_{\gamma} (B_{2^{-n+1} R}(x)) \le c_2  R^{\alpha -\delta} 2^{\alpha}    \sum_{n \ge 0} 2^{ n(\delta - \alpha)} < \infty. 
\end{eqnarray*}
Hence,
\[
|| R_1^G (1_{G} \cdot M^{\rho}_{\gamma}) ||_{L^{\infty} (E,m)} < \infty,
\]
and
\begin{equation}\label{eq;bdmro}
 \int_{G}  \int_{G} \frac{1}{|x-y|^{\delta}} \ M_{\gamma}^{\rho}(dy) M_{\gamma}^{\rho}(dx) < \infty. 
\end{equation}
By \cite[Exercise 4.2.2]{FOT}, \eqref{eq;bdmro} implies that $1_{G} \cdot M^{\rho}_{\gamma} \in S^G_{0}$.
Therefore, $1_{G} \cdot M^{\rho}_{\gamma} \in S^G_{00}$.
\qed

\begin{cor}
Almost surely in $X$, $M^{\rho}_{\gamma}$ does not charge capacity zero sets.
\end{cor}
\proof
Let $N \subset \R^2$ be an open set such that Cap$(N) =0$.
Note that by \cite[Lemma 2.2.3]{FOT}, $M^{\rho}_{\gamma}(1_{E_k} \cap N) = 0$  (see  Theorem \ref{t;s0ds}). Then the statement follows from 
\[ 
M^{\rho}_{\gamma}(N) \le \sum_{k \ge 1} M^{\rho}_{\gamma}(1_{E_k} \cap N) +  M^{\rho}_{\gamma}(\{0\}) = 0.
\]
\qed

The proof of the following theorem is a slight modification of \cite[Lemma 5.11]{ShTr13} in our setting.
\begin{thm}
Let $F_t^{k}$ be the positive continuous additive functional of $(X_t^{E_k})_{t \ge 0}$  in the strict sense associated to $1_{E_k} \cdot M_{\gamma}^{\rho} \in S_{00}^{E_{k}}$. Then, $F_t^{k} =  F_t^{k+1}, \; \forall t < \sigma_{E_{k}^{c}}$ $\P_{x}$-a.s. for all $x \in E_{k}$. In particular, $F_t : = \lim_{k \rightarrow \infty} F_t^{k}$, $t \ge 0$, is well defined in $A_{c,1}^{+,E}$, and related to $M_{\gamma}^{\rho}$ via the Revuz correspondence.
\end{thm}
\proof
We denote the set of all bounded, non-negative Borel measurable functions on $E_k$  by $\mathcal{B}_b^{+}(E_k)$.
Fix $f \in \mathcal{B}_{b}^{+}(E_{k})$ and for $x \in E_{k+1}$ define
\[
f_{k}(x) := \EE_{x} \Big[ \int_{0}^{\sigma_{E_k^{c}}} e^{-t} \, f(X_t) \,  d F_t^{k+1}  \Big].
\]
Since $f_{k} \in D(\E^{E_{k+1}})$ and $f_{k} =0$ $\E^{E_{k+1}}$-q.e. on $E_k^{c}$, we have $f_{k} \in D(\E^{E_{k}})$. For $x \in E_{k}$
\[
R_{1}^{E_k} \Big(f \ 1_{E_k} \cdot M_{\gamma}^{\rho}\Big)(x) = \EE_{x} \Big[ \int_{0}^{\sigma_{E_{k}^{c}}} e^{-t} \, f(X_t) \, d F_{t}^{k} \Big].
\]
Then, for $g \in \mathcal{B}_{b}^{+}(E_{k}) \cap L^2(E_{k}, m)$
\begin{eqnarray*}
\E^{E_{k}} \left( f_{k}, \, R_1^{E_k} g \right) & = & \E^{E_{k+1}} \left( f_{k}, \, R_1^{E_k} g \right) \\
& = &\int_{E}  R_1^{E_k}g \; f \,1_{E_{k+1}}  \cdot M_{\gamma}^{\rho} =
\E_1^{E_{k}} \left( R_{1}^{E_k} \Big(f  \,1_{E_k} \cdot M_{\gamma}^{\rho} \Big), \, R_1^{E_k}g  \right).\\
\end{eqnarray*}
Therefore, $f_{k} =  R_{1}^{E_k} \Big(f \,1_{E_k} \cdot M_{\gamma}^{\rho}\Big) \,$ $m$-a.e. Since $R_{1}^{E_k} \Big(f \,1_{E_k} \cdot M_{\gamma}^{\rho}\Big)$ is 1-excessive for $(R_{\alpha}^{E_k})_{\alpha >0}$, we obtain for any $x \in E_{k}$
\begin{eqnarray*}
R_{1}^{E_k} \Big(f \ 1_{E_k} \cdot M_{\gamma}^{\rho}\Big)(x) & = & \lim_{\alpha \rightarrow \infty} \alpha R_{\alpha+1}^{E_k} \left( R_{1}^{E_k} \Big(f \ 1_{E_k} \cdot M_{\gamma}^{\rho}\Big) \right)(x)\\
& = & \lim_{\alpha \rightarrow \infty} \alpha \int_{E_{k}} r_{\alpha+1}^{E_k}(x,y) \ R_{1}^{E_k} \Big(f \, 1_{E_k} \cdot M_{\gamma}^{\rho}\Big)(y) \ m(dy)\\
& = & \lim_{\alpha \rightarrow \infty} \alpha \int_{E_{k}} r_{\alpha+1}^{E_k}(x,y) \ f_{k}(y) \ m(dy)
= \lim_{\alpha \rightarrow \infty} \alpha R_{\alpha+1}^{E_k} f_{k}(x).
\end{eqnarray*}
Using in particular the strong Markov property, we obtain by direct calculation that the right hand limit equals $f_{k}(x)$ for any $x \in E_{k}$. Thus, we showed for all $x \in E_{k}$
\[
\EE_x \Big[ \int_{0}^{\sigma_{E_{k}^c}} e^{-t} \ f(X_t) \ d {F}_t^{k} \Big]
= \EE_x \Big[ \int_{0}^{\sigma_{E_{k}^c}} e^{-t} \ f(X_t) \ d {F}_t^{k+1} \Big].
\]
This implies that $F_t^{k} =  F_t^{k+1}, \; \forall t < \sigma_{E_{k}^{c}}$ $\P_{x}$-a.s. for all $x \in E_{k}$ (see e.g. \cite[IV. (2.12) Proposition]{BlGe}). Then, using Proposition \ref{p;cacod}, $F_t : = \lim_{k \rightarrow \infty} F_t^{k}$, $t \ge 0$, is well defined in $A_{c,1}^{+,E}$. Moreover, $(F_t)_{t\ge 0}$ is associated with $M_{\gamma}^{\rho}$ via the Revuz correspondence.
\qed 

Finally, almost surely in $X$, the time changed process  $(Z_t )_{t \ge 0}$ on $E$ can be defined as
\[
Z_t = X_{F_t^{-1}}, \quad t \ge 0,
\] 
where $F_t^{-1}:= \inf \{s>0 \mid F_s  > t  \}$, which is called the Liouville distorted Brownian motion.
\begin{remark}
In \cite{ShTr14} we considered the following assumptions: $\rho (x ) = |x|^{\alpha}$, $\alpha \in [2, \infty)$ and a symmetric (possibly) degenerate (uniformly weighted)  elliptic $2 \times 2$ matrix $A=(a_{ij})_{1 \le i,j \le 2}$, that is  $a_{ij} \in L^1_{loc} (\R^2,dx)$ and there exists a constant $\lambda \ge 1$
 such that for $dx$-a.e. $x \in \R^2$
\[
\lambda^{-1} \ \rho(x) \ \| \xi \|^2 \le \langle A(x) \xi, \xi \rangle \le \lambda \ \rho(x) \ \|\xi\|^2, \quad  \forall \xi \in \R^2,
\]
and the symmetric bilinear form
\[
\E^A(f,g) = \frac{1}{2} \int_{\R^2} \langle A \nabla f , \nabla g \rangle \, dx, \quad f,g \in C_0^{\infty}(\R^2).
\]

The closure $(\E^A,D(\E^A))$ of $(\E^A,C_0^{\infty}(\R^2))$ is a strongly local, regular, symmetric Dirichlet form. 
It is known from \cite{ShTr14} that there exists a  Hunt process $(Y_t)_{t \ge 0}$ starting from all points in $\R^2$ associated with the Dirichlet form $(\E^A,D(\E^A))$. Following the methods and techniques as in this section, we can construct the positive continuous additive functional  $(H_t)_{t\ge 0}$ of $(Y_t)_{t \ge 0}$ in the strict sense w.r.t. $M_{\gamma}^{\rho}$  and then the time changed process $Y_{H_t^{-1}}$, $t \ge 0$  on $E$ can be defined in the same way.
\end{remark}

\vspace*{2cm}
\noindent Jiyong Shin\\
School of Mathematics \\
Korea Institute for Advanced Study\\
85 Hoegiro Dongdaemun-gu,\\
Seoul 02445, South Korea,  \\
E-mail: yonshin2@kias.re.kr \\
\end{document}